\documentclass[a4paper,11pt]{article}
\usepackage{amssymb,latexsym}
\usepackage{amsmath}
\usepackage{amsfonts}
\usepackage{graphicx}
\usepackage{subfigure}
\usepackage{latexsym}
\usepackage{amsthm}
\usepackage{amssymb}

\theoremstyle{plain}

\newtheorem{theorem}{Theorem}

\theoremstyle{definition}

\newcommand{\R}{\mathbb{R}}
\newcommand{\inte}{{\rm int}}
\newcommand{\IR}{\mathbb R}

\newcommand{\F}{\mathcal F}

\newcommand{\be}{\begin{equation}}
\newcommand{\ee}{\end{equation}}
\newcommand{\ben}{\begin{enumerate}}
\newcommand{\een}{\end{enumerate}}

\begin{document}

\title{\textbf{Solving Inverse Problems for Steady-State Equations using A Multiple Criteria Model with Collage Distance, Entropy, and Sparsity}}
\author{Herb Kunze\footnote{Department of Mathematics and Statistics, University of Guelph, Guelph, Canada. Email: hkunze@uoguelph.ca} \and{Davide La Torre\footnote{SKEMA Business School - Universit$\acute{e}$ C$\hat{o}$te d'Azur, Sophia Antipolis Campus, France. Email: davide.latorre@skema.edu}} }
\maketitle

\begin{abstract}
In this paper, we extend the previous method for solving inverse problems for steady-state equations using the Generalized Collage Theorem by searching for an approximation that not only minimizes the collage error but also maximizes the entropy and minimize the sparsity. In this extended formulation, the parameter estimation minimization problem can be understood as a multiple criteria problem, with three different and conflicting criteria: The generalized collage error, the entropy associated with the unknown parameters, and the sparsity of the set of unknown parameters.
We implement a scalarization technique to reduce the multiple criteria program to a single criterion one, by combining
all objective functions with different trade-off weights. Numerical examples confirm that the collage method produces good, but sub-optimal, results. A relatively low-weighted entropy term allows for better approximations while the sparsity term decreases the complexity of the solution in terms of the number of elements in the basis.
\end{abstract}

\section{Introduction}
We present a multiple criteria model for solving an inverse problem for diffusion models described in terms of Partial Differential Equations (PDEs).
The theory of PDEs is crucial for modelling several problems in different disciplines, such as Business, Economics, Engineering, Finance and so on.
For instance it can be used to model the process of innovation and spread of ideas, the evolution of population dynamics, the dynamics of fluids, option pricing, and many other.

A PDE can be analyzed from both a direct and inverse approach: the direct problem is the analysis of the properties of existence, uniqueness, and stability of the solution. This is also refereed to the notion of well-posedness in the sense of Hadamard \cite{Hadamard}.

The inverse problem, instead, aims to identify causes from effects. In practice, this may be done by using observed data to estimate
parameters in the functional form of a model. Usually an inverse problem is ill-posed because some the properties related to existence, uniqueness, and stability
fail to hold.  When this happens, it is crucial to identify a suitable numerical scheme that ensures the convergence to the solution.

The literature is rich in papers studying ad hoc methods to address ill-posed inverse problems by minimizing a suitable approximation error along with utilizing some regularization techniques \cite{Kirsch,Neto,Tychonoff2,vogel}. The Collage-based approach instead, that has been utilized in this paper, relies
on an extension of the Collage Theorem \cite{Ba}, a consequence of Banach's fixed point theorem that has shown its importance to solve inverse problems for fixed point equations. The Collage Theorem is also the basis of the collage-based compression algorithm in fractal imaging \cite{Ba,KuLaMeVr}. It has also been extended
to inverse problems for ordinary differential equations and their application to different fields in \cite{KuVr99} and for partial differential equations
over solid and perforated domains in \cite{berkunlatrui,KuLaTo16}.

The results presented in this paper are a further contribution to this stream of research. We extend the Collage-based algorithm for steady-state equations by searching for an approximation that not only minimizes the collage error but also maximizes the entropy and minimize the sparsity. In this extended formulation, the parameter estimation minimization problem can be understood as a multiple criteria problem, with three different and conﬂicting criteria: The generalized collage error, the entropy associated with the unknown parameters, and the sparsity of the set of unknown parameters.
We implement a scalarization technique to reduce the multiple criteria program to a single criterion one, by combining all objective functions with different trade-off weights. Numerical examples confirm that the collage method produces good, but sub-optimal, results. A relatively low-weighted entropy term allows for better approximations while the sparsity term decreases the complexity of the solution in terms of the number of elements in the basis.

The paper is organized as follows: Section \ref{sec21} recalls some basic definitions in multiple criteria optimization.
Section \ref{sec2} recalls the main ideas of the Generalized Collage Theorem and how it can be used to solve
inverse problems for steady-state equations. Section \ref{sec3} presents the notion of entropy and how this formulation can be adapted to this particular context.
Section \ref{sec4} introduces the notion of sparsity and its importance to determine the complexity of the approximation.
Section \ref{sec6} formulates the multiple criteria model, which Section \ref{sec7} illustrates some numerical computations.
Section \ref{sec8} presents an application to an inverse problem in population dynamics and Section \ref{sec9}, as usual, concludes.

\section{Basics on Multiple Criteria Optimization}
\label{sec21}

This section focuses on recalling some basic facts in Multiple Criteria Optimization (MCO). In an abstract setting, a finite-dimensional MCO problem (see Sawaragi et al., 1985) can be stated as follows:
\begin{eqnarray}
\label{AP}
\max_{x\in X} \ J(x)
\end{eqnarray}
where $(X,\|\cdot\|)$ is a Banach space and $J:X\to\R^p$ is a vector-valued functional, and $\R^p$ is ordered by the Pareto cone $\R^p_+$. A point $x\in X$ is said to be Pareto optimal or efficient if $J(x)$ is one of the maximal elements of the set of achievable values $J(X)$. Thus a point $x$ is Pareto optimal if it is feasible and, for any possible $x'\in X$, $J(x)\le_{\R^p_+} J(x')$ implies $x=x'$. In a more synthetic way, a point $x\in X$ is said to be Pareto optimal if $(J(x) + \R^p_+) \cap J(X) = \{J(x)\}$.

\subsection{Scalarization}
Among the different techniques to reduce an MOP problem to a single criterion model there is, for sure, the scalarization technique.
Using a scalarization technique, a multiple objective model can be reduced to a single criterion problem by summing up all criteria with different weights.
The weights in front of each criterion express the relative importance of that criterion for the Decision Maker.
By using this approach, More precisely, by scalarization an MOP model boils down to:
\begin{eqnarray}
\label{SP}
\max_{x\in X} \sum_{i=1}^p \beta_i J_i(x)
\end{eqnarray}
where $\beta$ is a vector taking values in the interior of $\R^p_+$, namely $\beta\in\inte(\R^p_+)$. The equivalence between the scalarized problem and the original MOP problem is complete if the $J_i$ are linear and, by varying $\beta$, it is possible to obtain different Pareto optimal points.
In the other cases linear scalarization provides only partial results. Other scalarization methods can be found in the literature and one which is worth to be mentioned is the \textit{Chebyshev scalarization model} that can also be used for non-convex problems. Scalarization can also be applied to problems in which the ordering cone is different than the Pareto one. In this case, one has to rely on the elements of the dual cone to scalarize the multicriteria problem.

\subsection{$\epsilon$-constraint method}
The second model that is proposed to solve the vector-valued problem is the \textit{$\epsilon$-constraint method}.
In this methodology one of the objective functions is optimised using the others as constraints, then they are added to the constraint part of the model. The method is an hybrid methodology, in fact for the $\{J_i\}_{i\neq k}$, least acceptable levels, $\epsilon_i$ have to be set while the remaining objective function $J_k$ is optimised. Then the decision maker plays a relevant role in this model, choosing the objective function to be optimised and the least acceptable levels for the objective functions added as constraints.
Therefore, the original vector-valued problem can be now written as:

\begin{equation}
\max J_k(x)
\end{equation}
subject to:
\begin{equation}
\label{epsilon_system}
\left\{
\begin{array}{ll}
J_i(x) \geq \epsilon_i \ \ \ \ i\neq k\\
x\in X
\end{array}
\right.
\end{equation}

This method has the advantage of being theorically able to identify Pareto optimal points also of non-convex problems. However, it also has two potential drawbacks:
The identified optimal point is only granted to be weakly Pareto optimal, and the problem might become unfeasible due to the additional constraints.

\subsection{Goal Programming}
Another method that can be used to solve vector-valued problems that is worth to be mentioned is the \textit{Goal Programming} (or GP approach).
Goal Programming was first introduced by Charnes, Cooper, and Ferguson (1955) and Charnes and Cooper (1961). The innovative idea behind this model is the determination of the aspiration levels of an objective function. This model does not try to find an optimal solution but an acceptable one, it tries to achieve the goals set by the decision maker rather than maximising or minimising the objective functions.
Given a set of ideal goals \textit{$g_i$}, with \textit{$i=1, \dots,p$}, chosen by the decision maker, it is possible to re-write the problem in a GP form:
\begin{equation}
\min \sum_{i=1}^p \theta_i^+ \delta_i^+ + \theta_i^- \delta_i^- \nonumber
\end{equation}
Subject to:
\begin{equation}
\label{GPsystem}
\left\{
\begin{array}{l}
J_i(x)+\delta_i^- -\delta_i^+ = g_i   \ \ \ \  i=1,\dots, p\\
\delta_i^-,\delta_i^+ \geq 0\ \ \ \ \forall i=1,\dots,p\\
x\in X
\end{array}
\right.
\end{equation}

where $\delta_i^+$, $\delta_i^-$ are the positive and negative deviations (slack variables),respectively, and $\theta_i^+$,  $\theta_i^-$ are the corresponding weights.

\section{Inverse Problems for Steady-State Equations using the Generalized Collage Theorem}
\label{sec2}
The purpose of this paper is to provide an extended multiple criteria algorithm for the estimation of unknown parameters in steady-state equation by combining
the collage distance, the entropy, and the sparsity of the set of estimated coefficients.  The next Section \ref{sec8} will illustrate an application of this approach within the context of population dynamics.
Here we recall some basic facts about the collage-based approach. Before formulating the inverse problem for a generic family of problems in variational forms, let us have a quick look how a classical steady-state equation can be reformulated in such a form. Consider the steady-state equation: Find $u$ such that
$$
\left\{
  \begin{array}{ll}
  -\frac{d}{dx}\left(K(x)\frac{d u}{dx}(x)\right)= f(x), & x\in (0,1)\\
  u(0)=u(1)=0 &  \\
  \end{array}
\right.
$$
It is well-known that if we take the above model, we multiply both sides by a test function $\xi\in C^1_c([0,1])$ - the space of all continuous and differentiable functions with compact support in $[0,1]$ - and integrate over $[0,1]$,
the model can be written as
$$
\int_0^1 K(x)\frac{d u}{dx}(x) \frac{d \xi}{dx}(x)dx = \int_0^1 f(x) \xi(x) dx.
$$
If we define by $\lambda$ the vector of all unknown coefficients in $K$ and $f$, the model boils down to the more compact form:
$$
a_\lambda(u,\xi) = f^*_\lambda(\xi)
$$
where $a_\lambda$ is a family of bilinear forms defined as
$$
a_\lambda(u,\xi) = \int_0^1 K(x)\frac{d u}{dx}(x) \frac{d \xi}{dx}(x)dx
$$
and $f^*_\lambda$ is a family of linear operators defined as
$$
f^*_\lambda(\xi) = \int_0^1 f(x) \xi(x) dx
$$
The steady-state equation can be reformulated as a variational problem for which existence and uniqueness is guaranteed by the classical Lax-Milgram theorem
and over the space $H^1_0([0,1])$ of all integrable functions with integrable weak derivative on $[0,1]$ (see \cite{Evans} for more details in this).

\ \\
More in general, if $E$ is an Hilbert space, consider the following variational equation: Find $u\in E$ such that
\begin{equation}
a_\lambda(u,v) = f^*_\lambda(v),
\end{equation}
for any $v\in H$, where $f^*_\lambda(v)$ and $a_\lambda(u,v)$ are families of linear and bilinear maps, respectively, both defined on an Hilbert space $E$
for any $\lambda\in\Lambda$. Let $\langle\cdot\rangle$ denote the inner product in $E$, $\|u\|^2=\langle u,u\rangle$ and $d(u,v)=\|u-v\|$, for all $u,v\in E$.
The existence and uniqueness of solutions to this kind of equation are provided by the classical
Lax-Milgram representation theorem (see \cite{Evans}).
The following theorem presents how to determine the solution to the inverse problem for the above variational problem.
Following our earlier studies of inverse problems using fixed points of contraction mappings, we shall refer to it as a ``generalized collage method.''

\begin{theorem}(Generalized Collage Theorem)
\cite{KuLaToVr06-sub}
Let $\Lambda$ be a compact subset of $R^n$, and suppose that, for any $\lambda\in\Lambda$, $a_{\lambda}:E\times E\to \R$ be a family of bilinear forms and
$f^*_\lambda:E\to\IR$ be a family of linear forms. Furthermore, suppose that:

\begin{enumerate}
\item There exists a constant $M=\sup_{\lambda\in\Lambda} M_\lambda>0$ such that for any $\lambda\in\Lambda$, $|a_\lambda(u,v)|\le M_\lambda\|u\| \|v\|$ for all $u,v\in E$,
    \item There exists a constant $m=\inf_{\lambda\in\Lambda} m_\lambda>0$ such that for any $\lambda\in\Lambda$, $|a_\lambda(u,u)|\ge m_\lambda \|u\|^2$ for all $u\in E$.
\end{enumerate}

Then by the Lax-Milgram theorem, then for any $\lambda\in\Lambda$ there exists a unique vector
$u_{\lambda}$ such that
$$
a_{\lambda}(u_{\lambda},v) = f^*_\lambda(v)
$$
for all $v\in E$.  Then, for any $u\in E$,
\be
\|u-u_{\lambda}\|\le \frac{1}{m_{\lambda}}F(\lambda),
\label{collagedist2}
\ee
where
\be
\label{flambda}
F(\lambda)=\sup_{v\in E, ~ \|v\|=1}|a_{\lambda}(u,v)-f^*(v)| = \|a_\lambda(u,\cdot)-f^*\|.
\ee
\label{generalizedcollagetheorem}
\end{theorem}

In order to ensure that the approximation $u_{\mbox{\normalfont $\scriptstyle\lambda$}}$ is close to a target element
$u \in H$, we can, by the Generalized Collage Theorem, try to make the term $F({\mbox{\normalfont
$\lambda$}})/m_{\mbox{\normalfont $\scriptstyle\lambda$}}$ as close to zero as possible.
The appearance of the $m_{\mbox{\normalfont $\scriptstyle\lambda$}}$ factor complicates the procedure. However, if $\inf_{{\mbox{\normalfont $\scriptstyle\lambda$}}\in\Lambda}m_{\mbox{\normalfont
$\scriptstyle\lambda$}}\ge m>0$ then the inverse problem can be reduced to the minimization of the function
$F({\mbox{\normalfont $\lambda$}})$ on the space $\F$, that is,
\begin{equation}
\min_{{\mbox{\normalfont $\scriptstyle\lambda$}}\in\Lambda} F({\mbox{\normalfont $\lambda$}}).\label{minimizationproblem}
\end{equation}
In the following section we use the abbreviation CD, to denote the function $F({\mbox{\normalfont $\lambda$}})$.

\section{The Notion of Entropy}
\label{sec3}

The concept of entropy, as it is now used in information theory, was developed by C.E. Shannon \cite{Shannon}.
Over the years it has been used in different areas and applications in various scientific disciplines.
In his article, Shannon introduces the concept of information of a discrete random variable  with no memory as a functional that quantifies the uncertainty of a random variable. The concept of entropy describes the level of information associated with an event.
More precisely, the definition of Shannon's entropy \cite{Shannon,Flores} satisfies the following properties:

\begin{itemize}
  \item The measure is continuous and by changing the value of one of the probabilities by a very small amount should only produce a small change of the entropy;
  \item If all the outcomes are equally likely, then entropy should be maximal.
  \item If a certain outcome is a certainty, then the entropy should be zero.
  \item The amount of entropy should be the same independently of how the process is regarded as being
divided into parts.
\end{itemize}

According to these desiderata, Shannon defines the entropy in terms of a discrete random variable $X$, with possible outcomes $x_1,...,x_n$ as:
\be
ENT(X) = - \sum_{i=1}^n p(x_i) \ln(p(x_i))
\ee
For our purposes, this definition needs to be adapted to deal with a set of parameters, that can take both positive and negative values.
For a set of parameters ${\bf \lambda} = \{\lambda_1,\lambda_2,...,\lambda_n\}$ the notion of entropy is:
\be
ENT({\bf \lambda})=-\sum_{1}^{n} \frac{|\lambda_i|}{\lambda_T}\ln\frac{|\lambda_i|}{\lambda_T}
\ee
where $\lambda_T=\sum_i |\lambda_i|$. In the sequel, rather than maximizing the entropy term - that represents the total amount of information associated with
that particular combination of parameters' values - we will consider the minimization of its opposite, also known
as neg-entropy.
This criterion will be included in the multiple criteria model illustrated in the following Section \ref{sec6}.

\section{The Notion of Sparsity}
\label{sec4}
In literature the notion of sparsity has been widely used to reduce the complexity of a model by taking into in consideration only those parameters whose values have major impact on the solution. In other words, by adding this term we wish to determine solutions that are ``simple'', or more precisely \textit{sparse}.
We say that a real vector $x$ in $R^n$ is sparse, when most of the entries of $x$ vanish. We also say that a vector $x$ is $s$-sparse if it has at most $s$ nonzero entries. This is equivalent to say that the $\ell_0$-pseudonorm, or counting ‘norm’, defined as
\be
\|{\bf \lambda}\|_0=\#\{i:\lambda_i\neq 0\}
\ee
is at most $s$. The $\ell_0$-pseudonorm is a strict sparsity measure, and most optimization problems based on it
are combinatorial in nature, and hence in general NP-hard. To overcome these difficulties, it is common to replace the function with
relaxed variants or smooth approximations that measure and induce sparsity.
One possible variant is to use the $\ell_1$ norm instead, which is a convex surrogate for the $\ell_0$, defined as
\be
\|{\bf \lambda}\|_1=\sum_{i=1}^n |\lambda_i|
\ee
It is also the best surrogate in the sense that the $1$ ball is the smallest convex body containing all $1$-sparse objects of the form $±e_i$ (see \cite{Candes}).
Another possibility is to replace the $\ell_0$ pseudonorm with some approximation, as for instance
\be
\|\lambda\|_{*} = \sum_{i=1}^n \max\{e^{-\alpha \lambda_i},e^{\alpha \lambda_i}\}
\ee
or
\be
\|\lambda\|_{**} = \sum_{i=1}^n [\max\{e^{-\alpha \lambda_i},e^{\alpha \lambda_i}\}]^2
\ee
for a given $\alpha>0$. It is worth noticing that $\|\lambda\|_{**}$ is a $C^{1,1}$ or $LC^1$ function (continuous with Lipschitz gradient).

\section{The Model}
\label{sec6}
We now propose a new collage-based approach for solving inverse problems based on Multiple Criteria Optimization which combines together the Collage Distance,
the Entropy, and the Sparsity.  Then we consider the following criteria to be maximized/minimized simultaneously:

\begin{itemize}
  \item $CD(\lambda)$ is the Collage Distance, to be minimized over $\lambda\in\Lambda$. This criterion describes the accuracy of the approximation;
  \item $ENT(\lambda)$ is the Entropy, to be maximized over $\lambda\in\Lambda$. This criterion models the amount of information carried by the parameters' model;
  \item $SP(\lambda)$ is the Sparsity, to me minimized over $\lambda\in\Lambda$. This criterion instead describes the complexity of the solution in terms of number of elements in the basis to be utilized to approximate the target.
\end{itemize}

It is worth noticing that these three criteria are, in general, conflicting. It is clear that by reducing the sparsity criterion $SP(\lambda)$, this will
negatively affect the $CD(\lambda)$ as less elements in the basis are available to construct the solution.
To observe that the Entropy $ENT(\lambda)$ and the Sparsity $SP(\lambda)$ criteria are also conflicting, let us take a simple example where $X$ is a random variable with only two possible outcomes $x_1$ and $x_2$ with probabilities $p$ and $1-p$, respectively. It is clear that if $p$ increases, and then $1-p$ decreases, $x_1$ gets more and more likely to happen. This would produce a decrement in $ENT(X)$ while the sparsity of the vector $(x_1,x_2)$ would increase (see also \cite{Pastor} for a nice discussion on the importance of the concepts of entropy and sparsity).
By introduction the neg-entropy $-ENT$, the multiple criteria model can be formulated as a minimization program as follows:
\be
\min_{\lambda\in\Lambda} (CD(\lambda),-ENT(\lambda),SP(\lambda))
\ee
This multiple criteria problem can be transformed into a single criterion model by using one the approaches presented above.
In particular, one can construct the following single-criterion models:

\ \\
{\bf Model 1}: We scalarize the model by introducing three different positive weights, namely $\eta_1$, $\eta_2$, $\eta_3$.
The scalarized model boils down to:
$$
\min_{\lambda\in\Lambda} \eta_1 CD(\lambda) - \eta_2 ENT(\lambda) + \eta_3 SP(\lambda)
$$
The next section shows how the method works for different combinations of the weights.

\ \\
{\bf Model 2:} We move two out of three criteria into the constraints. Then Model 2 reads as:

$$
\min CD(\lambda)
$$
Subject to:
$$
\left\{\begin{array}{l}
  \lambda\in\Lambda \\
  -ENT(\lambda) \le \epsilon_1 \\
  SP(\lambda) \le \epsilon_2
\end{array}\right.
$$

\ \\
{\bf Model 3:} In the GP formulation, let us $g_1$, $g_2$, $g_3$ be the goals of $CD(\lambda)$, $ENT(\lambda)$, $SP(\lambda)$ respectively.
Then Model 3 reads as:
\begin{equation}
\min \sum_{i=1}^3 \theta_i^+ \delta_i^+ + \theta_i^- \delta_i^- \nonumber
\end{equation}
Subject to:
\begin{equation}
\label{GPsystem}
\left\{
\begin{array}{l}
  \lambda\in\Lambda \\
  CD(\lambda) + \delta_1^- -\delta_1^+ = g_1 \\
  -ENT(\lambda) + \delta_2^- -\delta_2^+ = g_2 \\
  SP(\lambda) + \delta_3^- -\delta_3^+ = g_3 \\
  \delta_i^-,\delta_i^+ \geq 0\ \ \ \ \forall i=1,\dots,3
\end{array}
\right.
\end{equation}

\section{A Computational Study}
\label{sec7}
To show the implementation of the algorithm, consider the following steady-state equation:
$$
\left\{
  \begin{array}{ll}
  -\frac{d}{dx}\left(K(x)\frac{d u}{dx}(x)\right)= f(x), & x\in (0,1)\\
  u(0)=u(1)=0 &  \\
  \end{array}
\right.
$$
with true $u(x)=x-x^2$, $d_K(x)=1+3x$, and $f(x)=-12x+1$.
The following tables show the result of our parameter estimation technique. Here we implement Model I presented in the previous section,
Models II and III can be implemented similarly.
Let us recall that $\eta_1$ is the coefficient of generalized collage distance, $\eta_2$ is the coefficient of entropy criterion, $\eta_3$ is the coefficient of the sparsity criterion.
In the following tables, we denote by $CD$ the value of the minimal general collage distance, by $ENT$ the value of the minimal entropy, and by $SP$ the value of minimal sparsity. $ER$ is the $L_2$ distance between the true $k(x)$ and the recovered $k(x)$.
We work with overlapping "hat" bases, the first with $11$ interior elements, the second with 23 interior elements, so every other element in the finer basis has the same peak point as an element from the coarser basis. We include the "half hats" at each end, as well, since we recover $K(x)$ in these bases and $K(x)$ is nonzero at the endpoints. So, in total there are $11+23+38$ elements.
In particular, the following Table 4 considers the effects of all three criteria simultaneously. This table shows how the three criteria interact differently when the three weights vary.

\begin{figure}
  \centering
  \includegraphics[width=5cm]{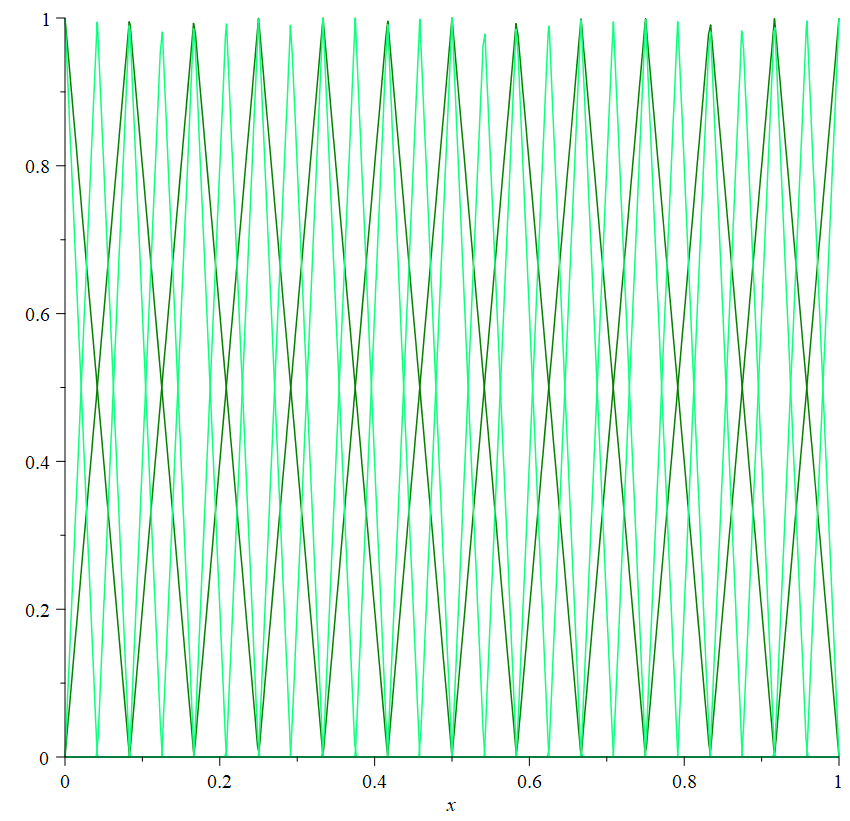}\\
  \caption{Hat basis with $11$ elements}
  \label{fig11}
\end{figure}

\begin{figure}
  \centering
  \includegraphics[width=5cm]{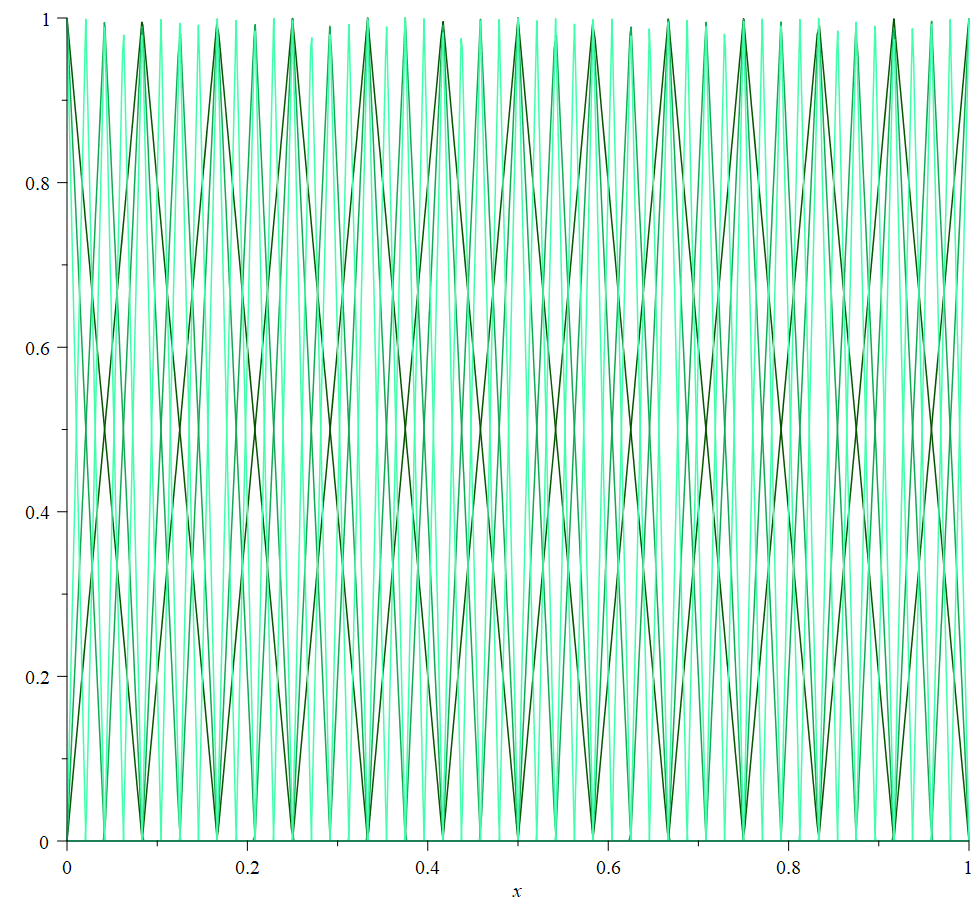}\\
  \caption{Hat basis with $23$ elements}
  \label{fig12}
\end{figure}

Table 1: CD versus ENT
$$
\begin{array}{cccccccc}
\eta_1    & \eta_2    & \eta_3    & CD                    & ENT       & SP & ER               \\
\hline
1.0 & 0.0 & 0 & 0.0000000000000000000 & -5.258363 & 38 & 0.12553297515980 \\
0.9 & 0.1 & 0 & 0.0000001334065773123 & -13.843577 & 38 & 0.04686792274056 \\
0.8 & 0.2 & 0 & 0.0000006753359361650 & -13.843580 & 38 & 0.04687624208382 \\
0.7 & 0.3 & 0 & 0.0000019845290201299 & -13.843584 & 38 & 0.04688698216096 \\
0.6 & 0.4 & 0 & 0.0000078272526691447 & -13.525527 & 38 & 0.20963394909766 \\
0.5 & 0.5 & 0 & 0.0000108023649479836 & -13.843597 & 38 & 0.04692168077792 \\
0.4 & 0.6 & 0 & 0.0000243008103020118 & -13.843607 & 38 & 0.04695245397844 \\
0.3 & 0.7 & 0 & 0.0000587837869956570 & -13.843625 & 38 & 0.04700459363829 \\
0.2 & 0.8 & 0 & 0.0001726457353481080 & -13.843661 & 38 & 0.04711204428741 \\
0.1 & 0.9 & 0 & 0.0008724077594118485 & -13.843769 & 38 & 0.04745936241855 \\
0.0 & 1.0 & 0 & 118740.8284896594930938881 & -13.979419 & 38 & 15.91754842300381 \\
\end{array}
$$

Table 2: CD vs sparsity
$$
{\small
\begin{array}{ccccccc}
\eta_1    & \eta_2    &     \eta_3    & CD                    & ENT       & SP & ER               \\
\hline
1.0 & 0 & 0.0 & 0.0000000000000000000 & -5.258363 & 38 & 0.12553297515980 \\
0.9 & 0 & 0.1 & 0.0005639498502964314 & -8.116850 & 33 & 2.70901013729197 \\
0.8 & 0 & 0.2 & 0.0000000000000000000 & -8.808666 & 31 & 3.19468535484908 \\
0.7 & 0 & 0.3 & 0.0000000000000000009 & -6.252997 & 26 & 1.03475376421398 \\
0.6 & 0 & 0.4 & 0.0042233007356215981 & -6.375100 & 32 & 1.53269905190502 \\
0.5 & 0 & 0.5 & 0.0000000000000000000 & -7.665059 & 29 & 2.71407610080125 \\
0.4 & 0 & 0.6 & 0.0000000000000000000 & -5.906960 & 32 & 0.00000000000000 \\
0.3 & 0 & 0.7 & 0.0065716738513130884 & -4.777920 & 28 & 0.98261846051097 \\
0.2 & 0 & 0.8 & 0.0000000000000000000 & -5.914810 & 27 & 1.94905796188513 \\
0.1 & 0 & 0.9 & 0.5203345122050530642 & -5.935058 & 30 & 2.53646236821924 \\
0.0 & 0 & 1.0 & 4111.1111111111111111111 & -0.000000 & 0 & 2.64575131106459 \\
\end{array}}
$$

Table 3: CD versus SP. Add another set of even finer basis elements, also sharing peak points with the other bases, for a total of 87 elements.
$$
\begin{array}{ccccccc}
\eta_1    &  \eta_2    &  \eta_3    & CD                    & ENT       & SP & ER               \\
\hline
1.0 & 0 & 0.0 & 0.0000000000000000000 & -8.432182 & 87 & 0.10275834304913 \\
0.9 & 0 & 0.1 & 0.0000000000000000000 & -9.285139 & 37 & 0.00000000000000 \\
0.8 & 0 & 0.2 & 0.0000000000000000000 & -8.582890 & 35 & 0.00000000000000 \\
0.7 & 0 & 0.3 & 0.0000000000000000000 & -8.052108 & 32 & 0.00000000000000 \\
0.6 & 0 & 0.4 & 0.0000000000000000000 & -7.115012 & 39 & 0.00000000000000 \\
0.5 & 0 & 0.5 & 0.0000000000000000000 & -4.397848 & 20 & 0.00000000000000 \\
0.4 & 0 & 0.6 & 0.0000000000000000000 & -8.002280 & 28 & 0.00000000000000 \\
0.3 & 0 & 0.7 & 0.0000000000000000000 & -7.709068 & 29 & 0.00000000000000 \\
0.2 & 0 & 0.8 & 0.0000000000000000000 & -5.204042 & 22 & 0.00000000000000 \\
0.1 & 0 & 0.9 & 0.0000000000000000000 & -7.699073 & 28 & 0.00000000000000 \\
0.0 & 0 & 1.0 & 4855.6857638888888888889 & -0.000000 & 0 & 2.64575131106459 \\
\end{array}
$$

Table 4: How the three criteria interact 
$$
\begin{array}{ccccccc}
\label{tablefinal}
\eta_1    &  \eta_2    &  \eta_3    & CD                    & ENT       & SP & ER               \\
\hline
0.96 & 0.02 & 0.02 & 0.0000005845708580369 & -30.713841 & 87 & 1.35681152916365 \\
0.38 & 0.02 & 0.60 & 0.0020363227800240392 & -20.407920 & 79 & 0.07009496798767 \\
0.58 & 0.02 & 0.40 & 0.0000007665555795281 & -15.963332 & 85 & 0.00008905869569 \\
\end{array}
$$

\begin{figure}
\centering
\begin{tabular}{cc}
  \includegraphics[width=4cm]{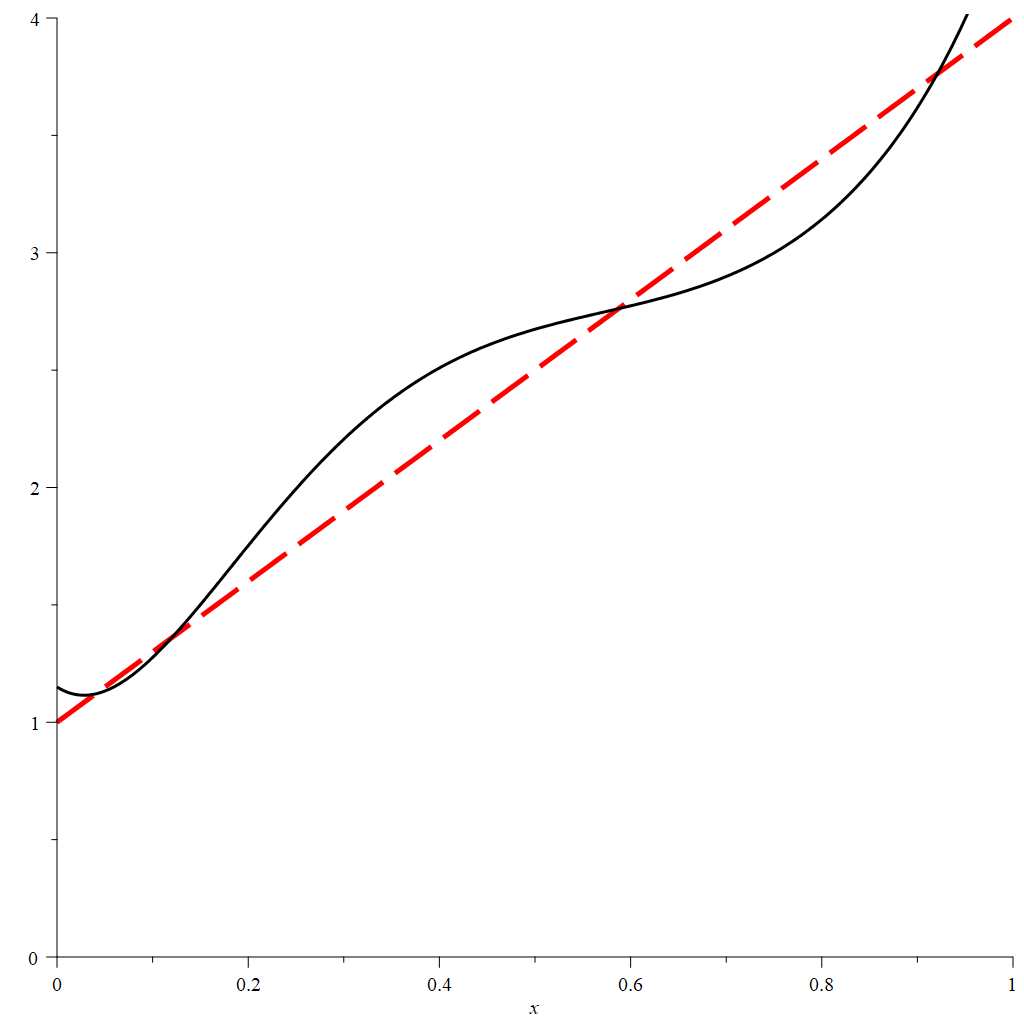} &   \includegraphics[width=4cm]{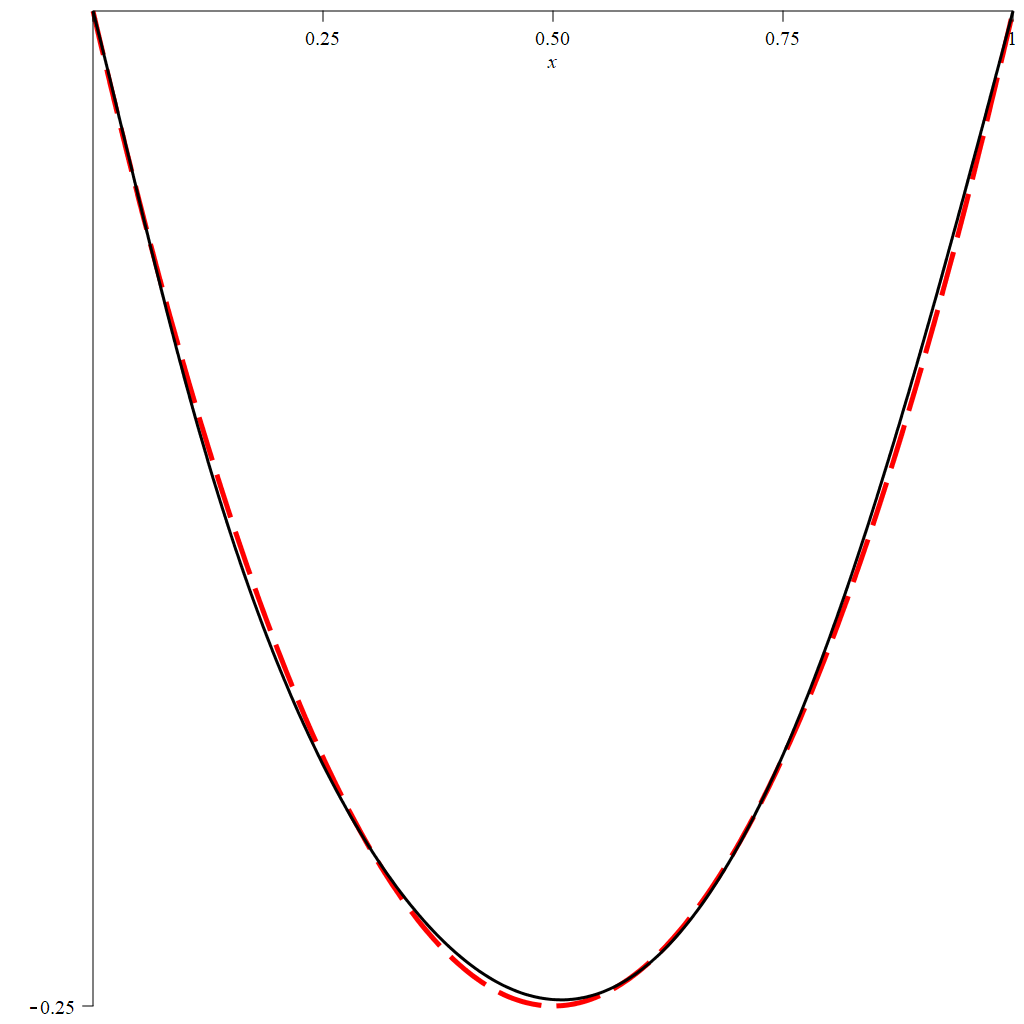} \\
\end{tabular}
\caption{True (K,u) vs Estimated (K,u),ent=0.02,sp=0.02}
\end{figure}

\begin{figure}
\centering
\begin{tabular}{cc}
  \includegraphics[width=4cm]{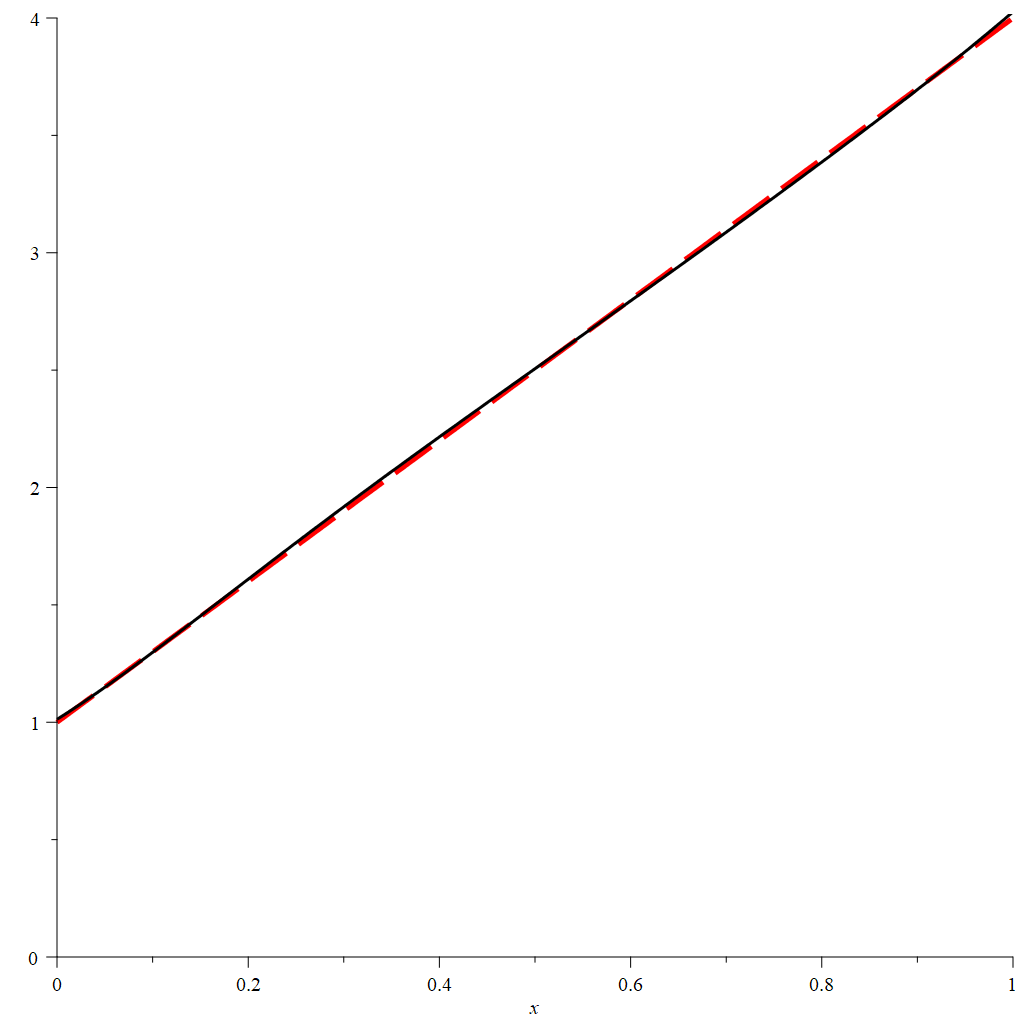} &   \includegraphics[width=4cm]{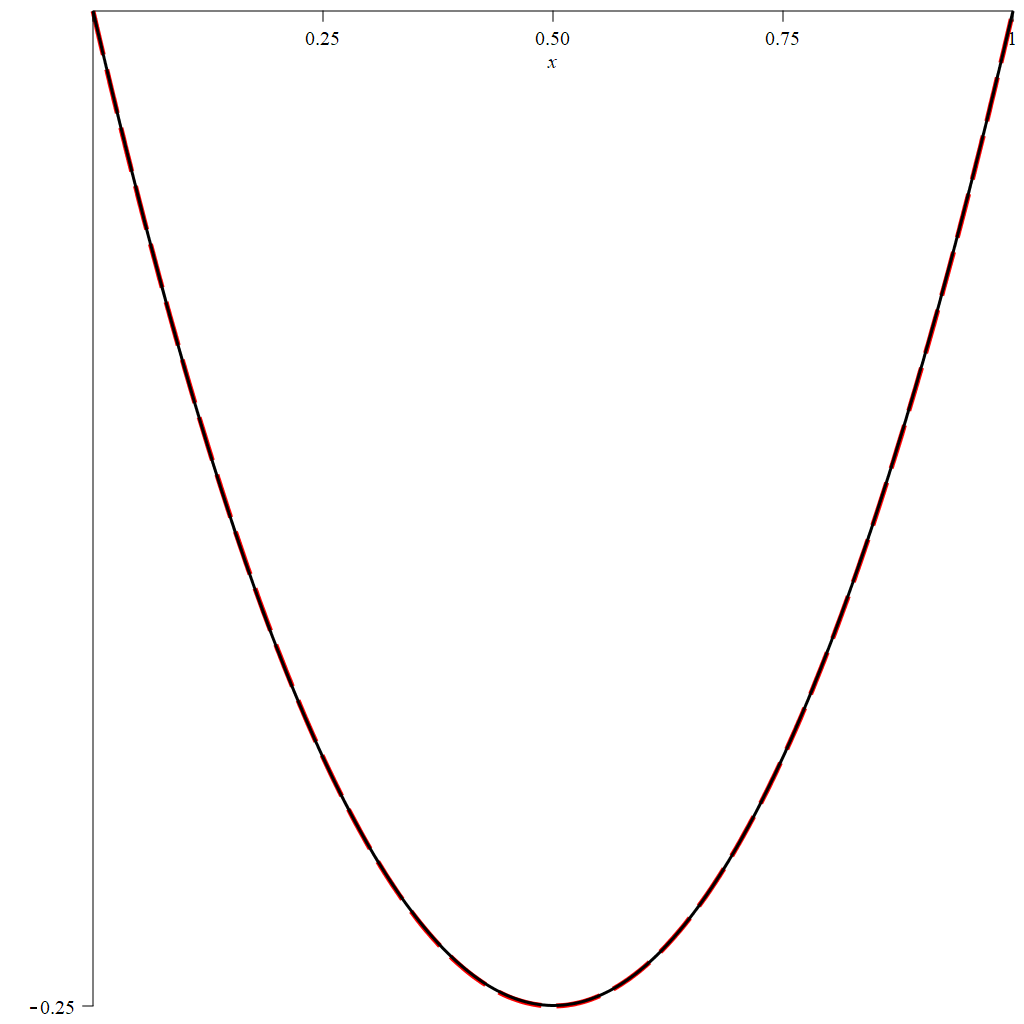} \\
\end{tabular}
\caption{True (K,u) vs Estimated (K,u),ent=0.02,sp=0.6}
\end{figure}

\begin{figure}
\centering
\begin{tabular}{cc}
  \includegraphics[width=4cm]{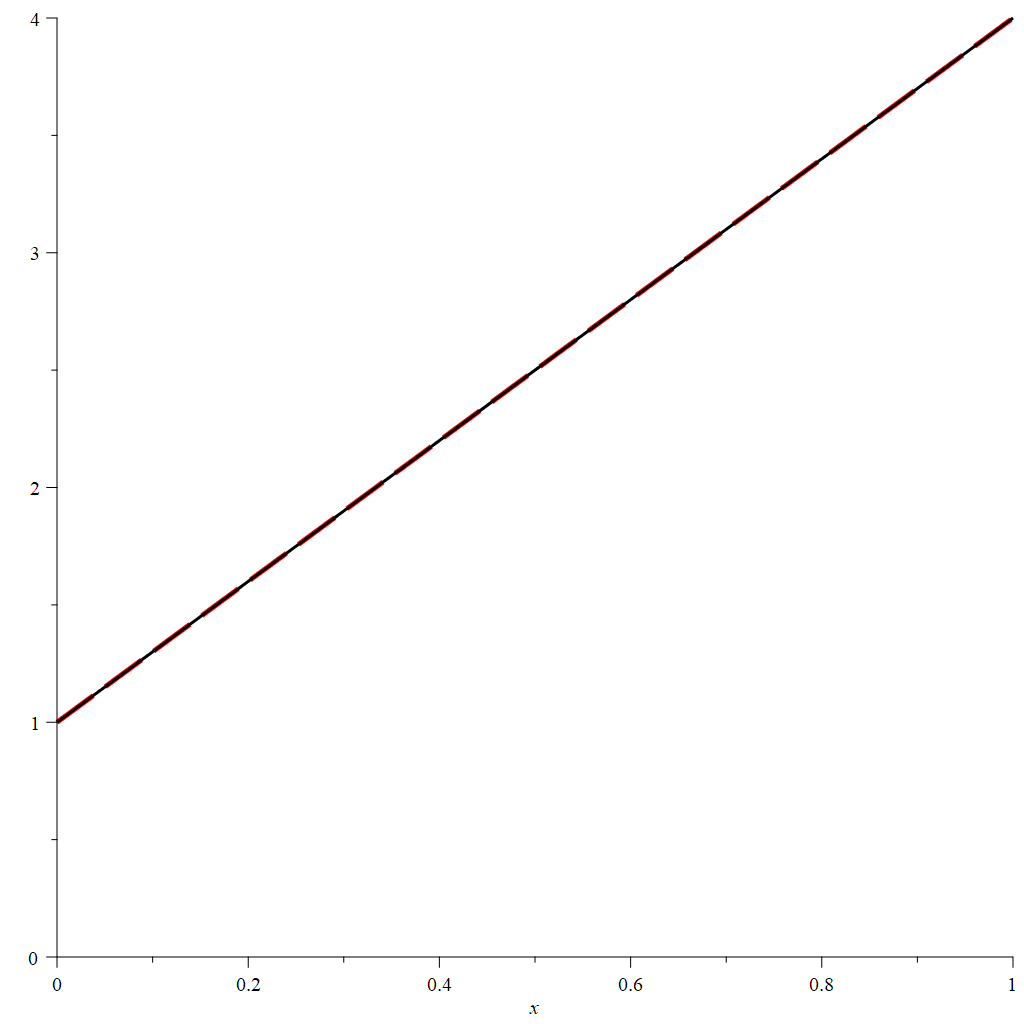} &   \includegraphics[width=4cm]{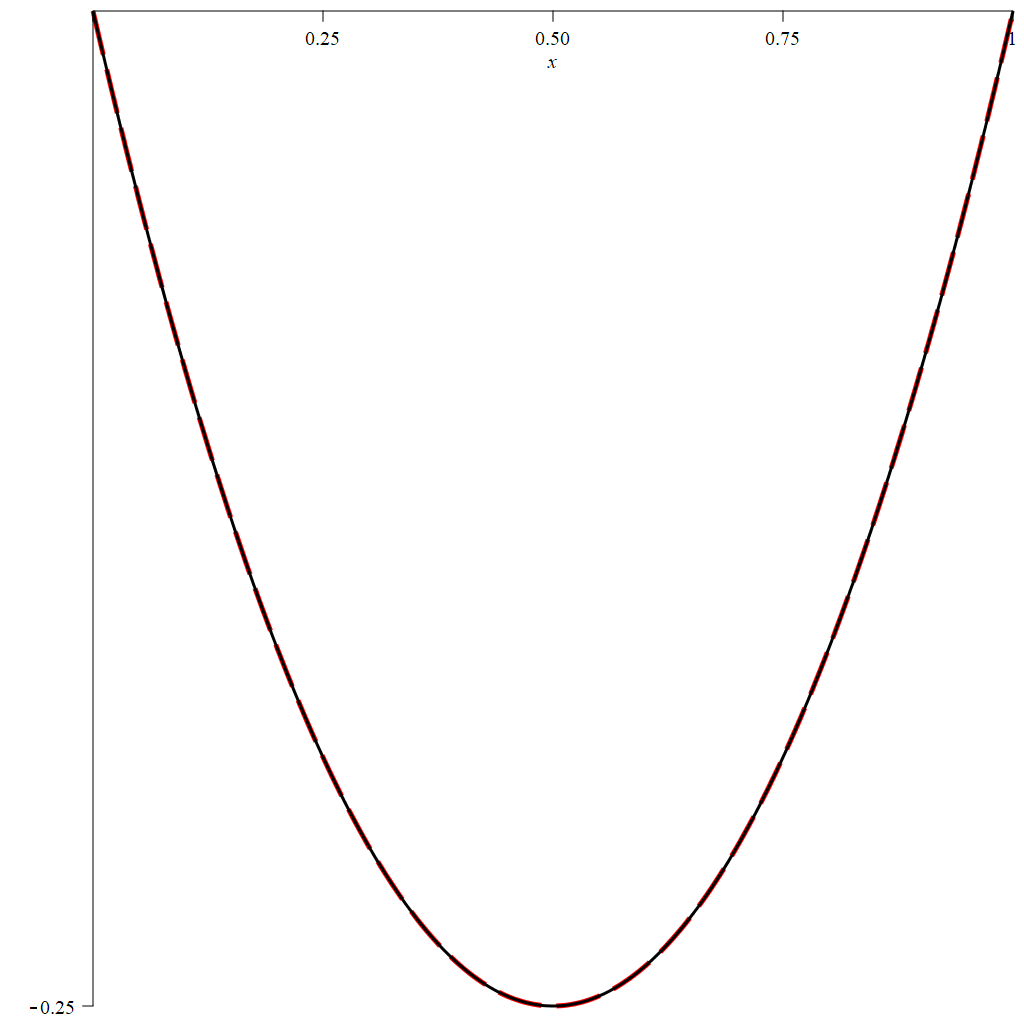} \\
\end{tabular}
\caption{True (K,u) vs Estimated (K,u),ent=0.02,sp=0.4}
\end{figure}

\section{An Application to Population Dynamics}
\label{sec8}
Population dynamics is an extremely important field and it is the basis of many economic models.
Recently a lot of attention has been devoted to a new stream of research, called Economic Geography, that aims at understanding
the evolution of population over space and time and the effects of migration flows on economies.
A spatial population model can be formulated as follows: Given a compact interval $[x_a,x_b]\subset R$, consider the following differential model with Dirichlet boundary conditions on $[x_a,x_b]$:
$$
\small{
\left\{
  \begin{array}{ll}
    {\partial P(x,t)\over \partial t} = \nabla \left(d_P(x) \nabla P(x,t) \right) + A(x), & (x,t)\in (x_a,x_b)\times (0,+\infty) \\
    P(x,t) = P_M, &  x\in \{x_a,x_b\}\\
    P(x,0) = P_0(x). & x\in [x_a,x_b]
  \end{array}
\right.}
$$
where $P(x,t)$ is the population level at time $t$ in location $x$, $P_0(x)$ is the initial distribution of population, and $A(x)$ is an exogenous flow of population. For simplicity, we also suppose that $P(x,t)$ is equal to a constant value $P_M$ at $x_a$ and $x_b$ and over time.
The steady-state level of population $\tilde P$ is the unique solution to the model
$$
\left\{
  \begin{array}{ll}
  -\nabla \left(d_P(x) \nabla \tilde P(x) \right) = A(x), & x\in (x_a,x_b)\\
  \tilde P(x) = P_M, &  x\in \{x_a,x_b\}\\
  \end{array}
\right.
$$

Let us consider the following numerical example where $x_a=0$, $x_b=1$, $d_P(x)=x+1$, $P_M=1$, and $A(x)=4x+1$.
In this case the true solution is $P(x)=x-x^2+1$.

\begin{figure}[h!]
	\begin{center}
        \includegraphics[width=9cm]{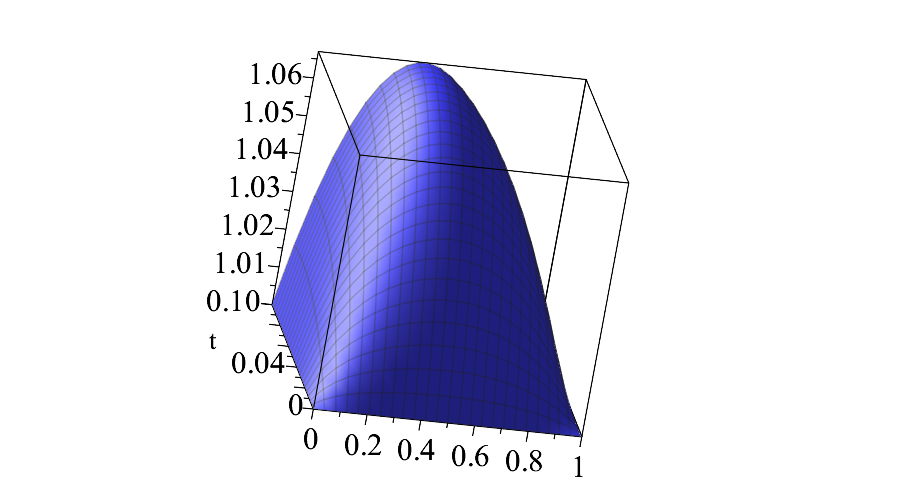}
		\caption{Evolution of $P$}
    \label{PFig}
	\end{center}
\end{figure}

The following tables present the results for inverse problem instead.
The first table shows the results with no noise added, $9$ interior data points, 11+23=34 interior basis functions.
The second table, instead, reports the result of the algorithm with $1\%$ relative noise added to data.

\begin{center}
\begin{tabular}{|c|c|c|c|c|c|c|c|}
\hline
$\eta_1$    & $\eta_2$    & $\eta_3$    & CD                    & ENT       & SP & ER  \\
\hline
1.0000000 & 0.0000000 & 0.0000000 & 0.00000000000000 & -3.817973 & 38 &  0.08124517 \\
0.9990000 & 0.0010000 & 0.0000000 & 0.00000000000351 & -13.935231 & 38 & 0.02629649 \\
0.9000000 & 0.1000000 & 0.0000000 & 0.00000004336279 & -13.935232 & 38 & 0.02630015 \\
0.8000000 & 0.2000000 & 0.0000000 & 0.00000021951333 & -13.935233 & 38 & 0.02630479 \\
0.4000000 & 0.6000000 & 0.0000000 & 0.00000789898350 & -13.935242 & 38 & 0.02634729 \\
0.1000000 & 0.9000000 & 0.0000000 & 0.00028361315324 & -13.935294 & 38 & 0.02663119 \\
0.9999000 & 0.0000000 & 0.0001000 & 0.00000000343665 & -4.083520 & 34 & 0.385780654 \\
0.9600000 & 0.0200000 & 0.0200000 & 0.00000000152452 & -13.935231 & 38 & 0.026297154 \\
0.7800000 & 0.0200000 & 0.2000000 & 0.00000224173370 & -6.525117 & 38 & 1.139612982 \\
0.5800000 & 0.0200000 & 0.4000000 & 0.00000282419676 & -5.958487 & 37 & 0.382813679 \\
\hline
\end{tabular}
\end{center}

\begin{center}
\begin{tabular}{|c|c|c|c|c|c|c|c|}
\hline
$\eta_1$    & $\eta_2$    & $\eta_3$    & CD                    & ENT       & SP & ER \\
\hline
1.0000000 & 0.0000000 & 0.0000000 & 0.00000000000000 & -3.925169 & 38 & 0.248001091 \\
0.9990000 & 0.0010000 & 0.0000000 & 0.00000000000422 & -13.934251 & 38 & 0.22401013 \\
0.9000000 & 0.1000000 & 0.0000000 & 0.00000005211046 & -13.934252 & 38 & 0.22400520 \\
0.8000000 & 0.2000000 & 0.0000000 & 0.00027783818395 & -4.462927 & 38 & 0.292032117 \\
0.4000000 & 0.6000000 & 0.0000000 & 0.00994433112542 & -4.473986 & 38 & 0.295136541 \\
0.1000000 & 0.9000000 & 0.0000000 & 0.34642831463774 & -4.538449 & 38 & 0.316445908 \\
0.9999000 & 0.0000000 & 0.0001000 & 0.00000000369833 & -3.272674 & 32 & 0.377684011 \\
0.9600000 & 0.0200000 & 0.0200000 & 0.00000000183208 & -13.934251 & 38 & 0.22400924 \\
0.7800000 & 0.0200000 & 0.2000000 & 0.00000144023205 & -5.710676 & 38 & 0.790857843 \\
0.5800000 & 0.0200000 & 0.4000000 & 0.00000012701174 & -11.825482 & 38 & 2.07316956 \\
\hline
\end{tabular}
\end{center}

\section{Conclusion}
\label{sec9}
The analysis of inverse problems for dynamical systems driven by differential equation is a crucial area in applied science.
In fact, in practical applications, it is relevant to be able to estimate the unknown parameters of a given equation starting from samples of the solution
collected by experiments or observations. In this paper we have extended a previous algorithm, based on the so-called "Collage Distance", to estimate
the unknown parameters of steady-state equations. This extended version also includes the notion of entropy and the notion of sparsity and it is modelled
as a multicriteria model. These three criteria are conflicting by nature, as a increment in precision usually implies an increment in sparsity.
We have solved the model using a scalarization technique with different weights and conducted several numerical experiments to show how the method works practically.


\begin{thebibliography}{9}

\bibitem{Ba}
Barnsley M, Fractals Everywhere, Academic Press, New York, 1989.

\bibitem{berkunlatrui} M.I. Berenguer, H. Kunze, D. La Torre, M. Ruiz Gal\'an, Galerkin method for constrained variational equations
and a collage-based approach to related inverse problems, J. Comput. Appl. Math. 292 (2016), 67--75.

\bibitem{Candes}
E. J. Candès (2014), Mathematics of sparsity (and a few other things),
Proceedings of the International Congress of Mathematicians, Seoul, South Korea, 2014.

\bibitem{Evans}
L.C. Evans (2010), Partial Differential Equations, Graduate Studies in Mathematics, American Mathematical Society.

\bibitem{Flores}
F. Flores Camachoa, N. Ulloa Lugob, H. Covarrubias Martıneza (2015), The concept of entropy, from its origins to teachers, Revista Mexicana de Fısica E 61, 69--80.

\bibitem{Kirsch}
A. Kirsch, An introduction to the mathematical theory of inverse problems, Springer, 2011.

\bibitem{KuVr99}
Kunze H and Vrscay E R, Solving inverse problems for ordinary differential equations using the Picard contraction mapping, Inverse Problems 15 (1999) 745--770.

\bibitem{kunze03a}
Kunze H and Gomes S, Solving An Inverse Problem for Urison-type Integral Equations Using Banach's Fixed Point Theorem, Inverse Problems 19 (2003) 411--418.

\bibitem{kunze03b}
Kunze H, Hicken J and Vrscay E R, Inverse Problems for ODEs Using Contraction Maps: Suboptimality of the ``Collage Method'', Inverse Problems 20 (2004)
977--991.




\bibitem{KuLaToVr06-sub}
Kunze H, La Torre D and Vrscay E R, A generalized collage method based upon the Lax--Milgram functional for solving boundary value inverse problems,  Nonlinear Anal. 71 12 (2009) e1337--e1343.


\bibitem{Ku3}
Kunze H, La Torre D, Vrscay E R, Solving inverse problems for DEs using the collage theorem and entropy maximization, Applied Mathematics Letters, 25 (2012), 2306-2311.

\bibitem{KuLaMeVr}
Kunze H, La Torre D, Mendivil F and Vrscay E R, Fractal-based methods in analysis, Springer, 2012.


\bibitem{KuLaTo16}
Kunze H and La Torre D, Collage-type approach to inverse problems for elliptic PDEs on perforated domains, Electronic Journal of Differential Equations, 48, 2015.





\bibitem{Hadamard}
J. Hadamard, Lectures on the Cauchy problem in linear partial differential equations, Yale University Press, 1923.


\bibitem{Neto}
F.D. Moura Neto, A.J. da Silva Neto, An Introduction to Inverse Problems with Applications, Springer, New York, 2013.

\bibitem{Pastor}
G. Pastor, I. Mora-Jimenez, R. Jantti, and A.J. Caamano (2013), Mathematics of Sparsity and Entropy: Axioms, Core Functions and Sparse Recovery,
Proceedings of the Tenth International Symposium in Wireless Communication Systems (ISWCS 2013).

\bibitem{Sawaragi} Sawaragi, Y., Nakayama, H., Tanino, T. (1985). Theory of multiobjective optimization (Academic Press, Inc.)



\bibitem{Shannon}
C.E. Shannon (1948), A Mathematical Theory of Communication, Bell System Technical Journal, 27 (3), 379–423.

\bibitem{Tychonoff2}
A.N. Tychonoff, N.Y. Arsenin, Solution of Ill-posed Problems, Washington: Winston \& Sons, 1977.

\bibitem{vogel}
C.R. Vogel, Computational Methods for Inverse Problems, SIAM, New York, 2002.

\end{thebibliography}
\end{document}